\documentclass{article}
\usepackage{amsmath}
\usepackage{amsfonts}
\usepackage{amssymb}	
\usepackage{color}
\usepackage{enumerate}
\usepackage{authblk}
\usepackage{graphicx}
\usepackage{comment}
\newtheorem{theo}{Theorem}[section]
\newtheorem{lem}[theo]{Lemma}
\newtheorem{cor}[theo]{Corollary}

\newtheorem{prop}[theo]{Proposition}
\newtheorem{defi}{Definition}[section]

\newcommand{\mysection}[1]{\section{#1} \setcounter{equation}{0}}
\newcommand{\proof}{{\sc Proof.} \quad}

\newcommand{\R}{\mathbb{R}}
\newcommand{\N}{\mathbb{N}}
\newcommand{\be}{\begin{equation} \label}
\newcommand{\ee}{\end{equation}}
\newcommand{\bes}{\begin{equation} \begin{array}{c} \label}
\newcommand{\ees}{\end{array} \end{equation}}
\newcommand{\bea}{\begin{eqnarray}\label}
\newcommand{\eea}{\end{eqnarray}}
\newcommand{\beas}{\begin{eqnarray} \begin{array}{rcl} \label}
\newcommand{\eeas}{\end{array} \end{eqnarray}}
\newcommand{\bas}{\begin{eqnarray*}}\newcommand{\eas}{\end{eqnarray*}}
\newcommand{\bass}{\begin{eqnarray*} \begin{array}{rcl}}
\newcommand{\eass}{\end{array} \end{eqnarray*}}
\newcommand{\basss}{\begin{eqnarray*} \begin{array}{c}}
\newcommand{\easss}{\end{array} \end{eqnarray*}}
\newcommand{\qed}{{}\hfill $\square$ \\}
\newcommand{\bit}{\begin{itemize}}
\newcommand{\eit}{\end{itemize}}
\newcommand{\nn}{\nonumber}
\newcommand{\eps}{\varepsilon}

\begin{document}
\title{Maximal dissipation in Hunter-Saxton equation for bounded energy initial data.}

\author[1]{Tomasz Cie\'slak}
\author[1,2]{Grzegorz Jamr\'oz}
\affil[1]{\small Institute of Mathematics, Polish Academy of Sciences, \'Sniadeckich 8, 00-656 Warszawa, Poland \newline{\small e-mail: T.Cieslak@impan.pl, G.Jamroz@impan.pl}}
\affil[2]{\small Institute of Applied Mathematics and Mechanics, University of Warsaw, Banacha 2, 02-097 Warszawa, Poland}

\maketitle

\begin{abstract} 
In \cite{zhang_zheng2005} it was conjectured by Zhang and Zheng that dissipative solutions of the Hunter-Saxton equation, which are known to be unique in the class of weak solutions,  dissipate the energy at the highest possible rate. 
The conjecture of Zhang and Zheng was proven in \cite{daf_jde} by Dafermos for monotone increasing initial data with bounded energy.  In this note we prove the conjecture in \cite{zhang_zheng2005} in full generality. To this end we examine the evolution of the energy of \emph{any} weak solution of the Hunter-Saxton equation. Our proof shows in fact that for every time $t>0$ the energy of the dissipative solution is not greater than the energy of any weak solution with the same initial data.  

\end{abstract}

\mysection{Introduction}\label{section1}

The Hunter-Saxton equation was introduced in \cite{hunter_saxton} as a simplified model to describe the evolution of perturbations of constant director fields in nematic liquid crystals. Liquid crystals share the mechanical properties of fluids and optical properties of crystals. Their description is essentially given by the evolution of two linearly independent vector fields, one denoting the fluid flow, another responsible for the dynamics of the so-called director fields, giving the orientation of the rod-like molecule. 

When analyzing planar director fields in the neighborhood of a constant field, after a nonlinear change of variables, one arrives at the following problem for $u:\R\times[0,T)\rightarrow \R$ 
\begin{equation}\label{hs}
\partial _tu(x,t)+\partial_x\left[\frac{1}{2}u^2(x,t)\right]=\frac{1}{2}\int_{-\infty}^xw^2(y,t)dy,\;\;w(x,t)=\partial_xu(x,t).
\end{equation}  
\[
u(x,0)=u_0(x), \;\;w(x,0)=w_0(x)=u_0'(x),\;\;-\infty<x<\infty. 
\]
Besides its physical meaning, equation \eqref{hs} possesses a few interesting mathematical properties. First, it seems to be a nice toy model for hydrodynamical equations. Next, it is completely integrable and admits infinitely many conservations laws, see \cite{hunter_zhengPhD}.

Due to its physical interpretation it is natural to expect from the solutions to \eqref{hs} that $\int_{\R}w^2(y,t)dy<\infty$, which means that the energy is finite. According to the Sobolev embedding in 1d, this means that for fixed time a solution to \eqref{hs} should be H\"{o}lder continuous in space, hence discontinuity does not occur. However, Lipschitz continuity might not be preserved. The well-known example, see \cite{bressan_cons},  \cite{daf_jhde}, of a solution admitting a blowup of the Lipschitz constant is given by 
\begin{equation}\label{z_P}
u(x,t) = \left\{\begin{array}{ccl}
0& 0\leq t<1  &-\infty<x\leq 0,\\ 
\frac{2x}{t-1}& 0\leq t<1 &0<x<(t-1)^2,\\
2(t-1)& 0\leq t<1 &(t-1)^2\leq x<\infty. \end{array}\right.
\end{equation}    
It develops a cusp singularity at $x=0, t=1$. However, it is possible to define global-in-time weak solutions. Such solutions have been constructed in \cite{hunter_zheng}. The question on the admissibility criteria yielding uniqueness appears. It was studied in \cite{hunter_zheng1}, \cite{zhang_zheng} and in the latter paper the notion of dissipative solutions was introduced. It was proven that these are unique solutions to \eqref{hs}. Such solutions were further studied in \cite{bressan_cons}. Let us now recall the definitions of both, weak and dissipative solutions.

\begin{defi}\label{weak}
A continuous function $u$ on $(-\infty,\infty)\times[0, \infty)$ is a weak solution of \eqref{hs} if, for each $t\in [0,\infty)$, $u(\cdot,t)$ is absolutely continuous on $\R$ satisfying
\[
\partial_xu(x,t)=w(x,t)\in L^\infty([0,\infty);L^2(\R)),
\] 
moreover the map $t\rightarrow u(\cdot,t)\in L^2_{loc}(\R)$ is absolutely continuous and \eqref{hs} holds on $\R\times[0, \infty)$, in the sense of distributions, $u(x,0)=u_0(x)$ and $w_0(x) \in L^2(\mathbb{R})$.
\end{defi} 
\begin{defi}\label{dissipative}
A weak solution to \eqref{hs} is called dissipative if its derivative $\partial_xu$ is bounded from above on any compact subset of the upper half-plane and $w(\cdot,t)\rightarrow w_0(x)$ in $L^2(\R)$ as $t\downarrow 0$. 
\end{defi}
As noticed in \cite{daf_jhde}, for any weak solution convergence of the solution to the initial data in $L^2$ is equivalent to the condition
\begin{equation}\label{en}
\limsup_{t\downarrow 0}E(t)\leq E(0),
\end{equation} 
$E(t):=\frac{1}{2}\int_{\R}w^2(y,t)dy$.
In \cite{daf_jde} Dafermos proved that dissipative solution is being selected by the criterion of maximal dissipation rate of the energy (or entropy, see \cite{daf_jde73}) among weak solutions for initial data $u_0$ being monotone increasing. He also stated that the same is true for general initial data with finite energy. In the present article we prove with all the details the above claim. It turns out that the proof in the general situation requires more involved reasoning. Our proof is based on the strategy of the proof in \cite{daf_jde}, however in order to handle the general situation we need an essentially more complicated argument. The biggest obstacle is that in order to proceed with the strategy of Dafermos, quite a detailed information on characteristics associated to any weak solution is required. Actually one needs to know how characteristics associated to weak solutions which are not disspative behave and how they push forward the energy. Our studies related to characteristics associated to weak solutions furnish enough information to enable us to execute the strategy. Most of the facts we prove on characteristics of weak solutions which are not dissipative seem to be new in the studies of \eqref{hs}. 

Let us now state the main results and recall an important definition from \cite{daf_jhde} which we shall need when dealing with the general case. 
\begin{theo}\label{Tw.1}
Let $u_0(x)$ be absolutely continuous with the derivative $u_0'(x)=w_0(x)$ a.e. and such that $w_0(x)\in L^2(\R)$. Then the dissipative solution of \eqref{hs} minimizes, for every $t \in [0,\infty)$, the energy among all weak solutions with the same initial data. This means that if $\tilde{u}$ is the unique dissipative solution of \eqref{hs} and $u$ any weak solution of \eqref{hs} starting from $u_0$, then
\begin{itemize}
\item $E^{\tilde{u}} (t) \le E^u (t)$ for every $t \ge 0$,
\item  if $E^{{u}}(t) = E^{\tilde{u}}(t)$ for every $t \ge 0$, then $u = \tilde{u}$.
\end{itemize}

\end{theo}
\begin{cor}

\label{Cor1}
The unique dissipative solution $\tilde{u}$ of \eqref{hs} maximizes the rate of the decay of energy among all weak solutions with the same initial data and consequently it is selected by the maximum dissipation principle. This means that if $u \neq \tilde{u}$ is a weak solution of \eqref{hs} such that $u(s) = \tilde{u}(s)$ for $s \in [0,t]$, then 
\begin{itemize}
\item there exists $s > t$ such that $E^u(s) > E^{\tilde{u}}(s)$ and
\item there is no $s>t$ such that $E^u(s) < E^{\tilde{u}}(s)$. 
\end{itemize}
\end{cor}

\begin{defi}\label{I_s}
For $s\in (0, \infty]$ we say that $I_s$ is a subset of the set 
\[
I=\{\zeta\in \R:u_0'(\zeta) \;\mbox{exists, equal to}\; w_0(\zeta)\}
\] 
consisting of such $\zeta\in I$ that $w_0(\zeta)>-\frac{2}{s}$. We denote $T_\zeta:=\infty$ if $w_0(\zeta)\geq 0$, $T_\zeta:=-\frac{2}{w_0(\zeta)}$ otherwise.  
\end{defi}

The paper is organized in the following way. In the next section we provide a sketch of the strategy of the proof of Theorem \ref{Tw.1}. The third section is devoted to introducing a collection of facts concerning characteristics. The fourth section is devoted to the study of energy contained between some pairs of characteristics. In the fourth section 
we study carefully the evolution of the positive part of the energy. Finally in the last section we formulate a proper averaging theory which enables us to arrive at the proof of Theorem \ref{Tw.1}. 

\mysection{The strategy of the proof of the main result}\label{section1.5}

In this section we describe with some details Dafermos' strategy of proving that dissipative solutions of \eqref{hs} are selected as the unique ones by the maximal energy dissipation criterion. It was successfully applied in the case of nondecreasing data in \cite{daf_jde}. We shall follow this strategy and very often we will be using some of the facts obtained in \cite{daf_jde}. Since exposing our result in a clear way requires a good source of reference concerning some of the computations done by Dafermos, we decided to recall many details of the latter in the present section. Finally, we will emphasize main additional difficulties which appear when one wants to execute the strategy in the case of absolutely continuous initial data with finite energy.

We notice that given a dissipative solution $u$, see Definition \ref{dissipative}, by \cite[Theorem 4.1]{daf_jhde}, we know that 
\begin{equation}\label{4.1}
\int_{\R}w^2(y,t)dy=\int_{I_t}w_0^2(\zeta) d\zeta.
\end{equation}
Notice that if we prove that the energy $E=\frac{1}{2}\int_{\R}w^2(y,t)dy$ associated to any weak solution of \eqref{hs} is bounded from below by  $\int_{I_t}w_0^2(\zeta) d\zeta$ and any weak solution $u$ satisfying \eqref{4.1} is a disspative solution, then we are done. Hence, in order to complete the proof of Theorem \ref{Tw.1} it is enough to prove the following two propositions.
\begin{prop}\label{prop2.1}
Let $u$ be a weak solution of \eqref{hs} (see Definition \ref{weak}). Moreover, assume \eqref{4.1} holds. Then $u$ is actually a dissipative solution.
\end{prop}
\begin{prop}\label{prop2.2}
Let $u$ be a weak solution of \eqref{hs}. Then
\begin{equation}\label{ostatnie}
\int_{\R}w^2(y,t)dy\geq\int_{I_t}w_0^2(\zeta) d\zeta.
\end{equation}
\end{prop}  
Now, we shall recall how the strategy outlined above was executed by Dafermos for nondecreasing initial data. Thus, we will be able to explain to the reader what difficulties appear for more general initial data. Moreover, some formulas which appear in this section will be used by us in a more complicated framework, and it seems to us useful to introduce them in the basic setting.

A characteristic associated to the weak solution $u$ of \eqref{hs} is a Lipschitz continuous function $x:[0,T]\rightarrow \R$ satisfying 
\begin{equation}\label{2.0.2}
\dot{x}(t)=u(x(t),t)\;\mbox{for a.e.}\; t\in[0,T],
\end{equation}
\[
x(0)=x_0.
\]
By \cite[Lemma 3.1]{daf_jhde} we know that for every $x_0$ there exists a characteristic $x(t)$ of \eqref{hs} (perhaps not unique) passing through $x_0$. Morever, every characteristic is actually a $C^1$ function and satisfies
\begin{equation}\label{2.0.01}
\dot{x}=u(x(t),t):=u_{x(t)}(t),\;\;\; \dot{u}_{x(t)}(t)=\frac{1}{2}\int_{-\infty}^{x(t)}w^2(y,t)dy
\end{equation}   
pointwise and a.e., respectively. The function $t\rightarrow u(x(t),t)$ is Lipschitz continuous. 

Following \cite{daf_jde}, given characteristics $x_1(t), x_2(t)$ emanating from $x_1, x_2\in \R$ we introduce 
\begin{eqnarray}\label{hpt}
h(t):=x_2(t)-x_1(t)&,& p(t):=u(x_2(t),t)-u(x_1(t),t), \\
\omega^{x_1, x_2}(t)&:=&\frac{p(t)}{h(t)}.
\end{eqnarray}
One sees that 
\begin{equation}\label{hpt1}
\dot{h}(t)=p(t), \; \dot{p}(t)=\frac{1}{2}\int_{x_1(t)}^{x_2(t)}w^2(y,t)dy.
\end{equation}   
An immediate consequence of \eqref{hpt1} is that if $p$ is initially positive, then $h$ stays positive during its evolution. In other words, nondecreasing initial data assure that characteristics do not intersect. Since $h(t)>0$ for $t\geq 0$,
\[
\dot{h}(t)=\omega^{x_1, x_2}(t)h(t)
\]
and so
\begin{equation}\label{wazne}
h(t)=h(0)e^{\int_0^t\omega^{x_1, x_2}(s)ds}.
\end{equation}
Moreover, 
\begin{equation}
\label {omegaintegral}
\dot{\omega}^{x_1, x_2}(t)=-\left(\omega^{x_1, x_2}\right)^2(t)+\frac{1}{2h(t)}\int_{x_1(t)}^{x_2(t)}w^2(y,t)dy.
\end{equation}
Next, Dafermos shows that 
\begin{equation}\label{hpt2}
\int_{x_1(t)}^{x_2(t)}w^2(y,t)dy\geq \int_{x_1}^{x_2}(w_0)^2(y)dy
\end{equation}
and in view of the fact that $h(t)>0$, which implies that $\R\setminus I_t$ is of measure $0$, concludes the proof of Proposition \ref{prop2.2}. On the other hand if $E(t)=E(0)$ for any $t>0$, then \eqref{hpt2} must also hold as equality for any pair of characteristics $x_1(t), x_2(t)$. This leads Dafermos to the fact that $u$ must be a dissipative solution, see \cite[the end of section 3]{daf_jde}. So Proposition \ref{prop2.1} also holds, thus Theorem \ref{Tw.1} is true for nondecreasing initial data. 

Now, let us comment on the difficulties which appear when considering general initial data. First of all, \eqref{hpt1} does not guarantee that characteristics do not intersect. Actually, the collision of characteristics is possible. Next, it is also possible that characteristics of weak solutions branch. Moreover, our proof requires treating separately the positive and negative parts of the energy as well as change of variables formulas for Lebesgue-Stieltjes integral. It is enough to notice that a solution given in \eqref{z_P} can be continued for the times $t>1$ in the following non-unique way.
\begin{equation}\label{z_Q}
u(x,t) = \left\{\begin{array}{ccl}
0& t>1  &-\infty<x\leq 0,\\ 
\frac{2x}{t-1}& t>1 &0<x<k(t-1)^2,\\
2k(t-1)& t>1 &k(t-1)^2\leq x<\infty, \end{array}\right.
\end{equation}
$k\geq 0$. In order to deal with those obstacles and execute the strategy of Dafermos in the case of absolutely continuous initial data, we need detailed studies of characteristics of weak solutions which may collide and branch, which is done in the next section. 
Finally, we observe that both Proposition \ref{prop2.1} and Proposition \ref{prop2.2} are consequences of the following one.
\begin{prop}\label{prop2.3}
Let $u$ be a weak solution of \eqref{hs}. Let $x_\xi$ and $x_\zeta$ be characteristics starting at $\xi, \zeta \in \R$, respectively, $\xi<\zeta$. Moreover, assume $x_\xi(t)\leq x_\zeta(t)$ for any $t\geq 0$. Then the following formula holds
\begin{equation}\label{ostatnie_po_char}
\int_{x_\xi(t)}^{x_\zeta(t)}w^2(y,t)dy\geq\int_{I_t\cap (\xi, \zeta)}w_0^2(\zeta) d\zeta\;.
\end{equation}
\end{prop}
Proposition \ref{prop2.2} is implied by Proposition \ref{prop2.3} in a straightforward way. To see the implication from 
Proposition \ref{prop2.3} to Proposition \ref{prop2.1} one needs to take into account that, as was stated in Section \ref{section2}, to any weak solution a set of $C^1$ characteristics is associated. They however might not be unique. Clearly,
for any characteristic $x_\xi(t)$ emanating from a point $\xi$, we have (see \eqref{2.0.01})
\begin{equation}\label{raz}
\dot{u_\xi}(t)=\frac{1}{2}\int_{-\infty}^{x_\xi(t)}w^2(y,t)dy,
\end{equation}
where $u_\xi(t):=u(x_\xi(t),t)$. 

On the other hand, we see that any weak solution to \eqref{hs} satisfying \eqref{4.1} and \eqref{ostatnie_po_char}, satisfies also
\begin{equation}\label{dwa}
\int_{x_\xi(t)}^{x_\zeta(t)}w^2(y,t)dy=\int_{I_t\cap (\xi, \zeta)}w_0^2(\zeta) d\zeta.
\end{equation} 
Hence, taking into account \eqref{dwa} and integrating \eqref{raz} in time
we arrive at
\[
u(x_\xi(t),t)=u_0(\xi)+\frac{1}{2}\int_0^t\int_{I_s\cap(-\infty,\xi)}w_0^2(\zeta)d\zeta ds.
\]
The above equality tells us that $u$ is actually a dissipative solution to \eqref{hs} according to \cite[Theorem 2.1]{bressan_cons}, see also \cite[Theorem 2.1]{daf_jhde}. 

In view of the above, all we have to show, in order to complete the proof of Theorem \ref{Tw.1}, is Proposition \ref{prop2.3}.

\mysection{Some information on characteristics}\label{section2}

In this section we study the behavior of characteristics associated to weak solutions of \eqref{hs}. First we notice the following lemma.
\begin{lem}\label{lemat1}
Consider any weak solution $u$ of \eqref{hs}. Let $\xi_0 \in I_t$. Choose $\xi_1\in \R$ such that $|\xi_1-\xi_0|$ is small enough. Then for any $x_{\xi_0}(t), x_{\xi_1}(t)$, characteristics associated to $u$ emanating from $\xi_0, \xi_1$ respectively, there exists a positive continuous function $\chi$ such that
\begin{equation}\label{trzy}
|x_{\xi_1}(s)-x_{\xi_0}(s)|\geq \chi(s)\;\;\mbox{for}\;s\in(0,t].
\end{equation}
\end{lem}
\proof
We take two characteristics emanating from $\xi_0$ and $\xi_1$. In view of \eqref{wazne},
as long as $\int_0^t\omega(s)ds>-\infty$, $h(t)>0$ and so characteristics do not intersect.

Moreover, \eqref{omegaintegral} is satisfied
as long as $h>0$. It implies, by the Schwarz inequality, 
\begin{equation}\label{omega}
\dot{\omega}^{x_{\xi_0}, x_{\xi_1}}(t)\geq -\frac{1}{2}\left(\omega^{x_{\xi_0}, x_{\xi_1}}(t)\right)^2, \;\;\mbox{hence}\;\;\omega^{x_{\xi_0}, x_{\xi_1}}(t)\geq \frac{2\omega^{\xi_0, \xi_1}(0)}{2+t\omega^{\xi_0, \xi_1}(0)}.
\end{equation}
Next, we take $\xi_0\in I_t$ and observe that for $\xi_1$ in a sufficiently close neighborhood of $\xi_0$,
\[
\omega^{\xi_0, \xi_1}(0)>-\frac{2}{t}\;.
\]
Indeed, since $\xi_0\in I_t$ we have on the one hand $w(\xi_0)>-\frac{2}{t}+\varepsilon_0$ for some small $\varepsilon_0$, and on the other hand,
\[
\frac{u(x)-u(\xi_0)}{x-\xi_0}=u'(\xi_0)+\frac{o(|x-\xi_0|)}{x-\xi_0}\;.
\] 
We choose $\xi_1$ such that $\frac{o(|\xi_1-\xi_0|)}{\xi_1-\xi_0}<\frac{\varepsilon_0}{2}$.   Then, 
\[ \omega ^{{\xi _0, \xi_1}}(0) > -\frac 2  t + \frac {\eps_0}{2}. \]
Hence, in view of \eqref{omega}, \eqref{wazne} leads to

\begin{eqnarray*}
|x_{\xi_1}(s) - x_{\xi_0}(s)| &=& |\xi_1 - \xi_0| \exp\left({\int_0^s \omega^{x_{\xi_0}, x_{\xi_1}}(\tau)d\tau}\right)\\ &\ge& |\xi_1 - \xi_0| \exp \left({\int_0^s \frac{2\omega^{\xi_0, \xi_1}(0)}{2+\tau\omega^{\xi_0, \xi_1}(0)}d\tau }\right) \\&=&  |\xi_1 - \xi_0| \frac 1 4 (2+s \omega^{\xi_0,\xi_1}(0))^2 \\&\ge& |\xi_1 - \xi_0| \frac 1 4 (2+t \omega^{\xi_0,\xi_1}(0))^2 \ge \frac 1 {16}|\xi_1 - \xi_0|  t^2 \eps_0^2.
\end{eqnarray*}


\qed

As a corollary we infer the following fact.
\begin{cor}\label{cor2.1}
Consider a characteristic $x(t)$, associated to a weak solution $u$ of \eqref{hs}, emanating from $x_0\in I_t$. This characteristic does
not cross with any other until time $t$.
\end{cor}
\proof
Any characteristic starting from the neighborhood of $x_0$ does not cross $x(t)$ by Lemma \ref{lemat1}. Next consider a characteristic starting from a point being outside of a neighborhood of $x_0$. If it crosses $x(t)$ then, in particular, it crosses one of the characteristics from the neighborhood of $x_0$. But this way we obtain a characteristic starting from a neighborhood of $x_0$ which intersects $x(t)$, which leads to contradiction.  

\qed

However, as we have seen in \eqref{z_Q}, characteristics associated to a weak solution can branch. We need to find out how often it may happen in order to proceed with the proof. We have the following lemma.
\begin{lem}\label{lemat2}
Let $u$ be a weak solution of \eqref{hs}. For almost every (more precisely, all except a countable number) $x_0\in I_T$, the
characteristic associated to $u$ starting from $x_0$, does not branch before time $T$.
\end{lem} 
Combining Lemma \ref{lemat2} with Lemma \ref{lemat1} we obtain the following claim.
\begin{cor}\label{cor2.2}
Let $u$ be a weak solution to \eqref{hs}. For every except a countable number $x_0\in I_t$, the characterictic emanating from $x_0$ is unique forwards and backwards up to time $t>0$.
\end{cor}
The proof of Lemma \ref{lemat2} requires some steps, in particular the introduction of leftmost and rightmost characteristics. To this end we show a few claims.   
\begin{prop}\label{prop2.0.1}
Consider a family of Lipschitz continuous functions $f_\alpha:[a,b]\rightarrow \R$, $\alpha\in A\subset \R$ satisfying
\begin{equation}\label{2.0.1}
|f_\alpha(t)-f_\alpha(s)|\leq L|t-s|
\end{equation}
for all $t,s\in [a,b]$ and some $L>0$. Then both $\sup_{\alpha\in A}f_\alpha$ and $\inf_{\alpha\in A}f_\alpha$ are Lipschitz continuous.
\end{prop}
\proof
We shall prove the claim of the proposition only for $\sup f_\alpha$, for $\inf f_\alpha$ the proof is the same.
First we fix $t,s\in[a,b]$, $t>s$. We notice that for any $\varepsilon>0$ there exist $\alpha_0, \alpha_1\in A$ such that  
\[
\sup_{\alpha\in A}f_\alpha(t)-\varepsilon< f_{\alpha_0}(t),
\]
\[
\sup_{\alpha\in A}f_\alpha(s)-\varepsilon< f_{\alpha_1}(s).
\]
Hence
\[
f_{\alpha_1}(t)-f_{\alpha_1}(s)< \sup_{\alpha\in A}f_\alpha(t)-\sup_{\alpha\in A}f_\alpha(s)+\varepsilon
\]
and
\[
f_{\alpha_0}(t)-f_{\alpha_0}(s)> \sup_{\alpha\in A}f_\alpha(t)-\sup_{\alpha\in A}f_\alpha(s)-\varepsilon,
\]
which together with \eqref{2.0.1} yields
\[
-L(t-s)-\varepsilon\leq f_{\alpha_1}(t)-f_{\alpha_1}(s)-\varepsilon< \sup_{\alpha\in A}f_\alpha(t)-\sup_{\alpha\in A}f_\alpha(s)< f_{\alpha_0}(t)-f_{\alpha_0}(s)+\varepsilon\leq L(t-s)+\varepsilon.
\]
Letting $\varepsilon$ go to $0$ in the above inequality we obtain the claim of the proposition.

\qed

Let us now state a proposition, which we will use in the sequel, which is a consequence of the Kneser's theorem, see \cite[Theorem II.4.1]{hartman}, as well as the fact that any weak solution to \eqref{hs} is continuous. 
\begin{prop}\label{kneser}
Let the image of a point under characteristics emanating from $(x_0,t_0)$ be defined as
\[
A(t):=\{(z,t):z=x(t), x(t_0)=x_0, x(t)\;\mbox{is a characteristic of}\;\eqref{hs}\;\mbox{associated to the weak solution}\;u \}.
\]
Then $A(t)$ is a compact and connected set.
\end{prop}
The next lemma contains the proof of existence of rightmost and leftmost characteristics.
\begin{lem}\label{lemat3}
Let $u : [0, T]\times \R\rightarrow \R$ be a bounded continuous function solving \eqref{hs} in the weak sense and let $x_\alpha(t)$ be a family of characteristics associated to $u$, i.e.
$C^1$ functions on $[0,T]$ satisfying \eqref{2.0.2}.
Then function $y : [0,T]\times \R\rightarrow \R$ defined for $t\in [0,T]$ by $y(t):=\sup_{\alpha}x_\alpha(t)$ is also a characteristic of \eqref{hs} associated to the weak solution $u$. The same claim holds for $\inf_{\alpha}x_\alpha(t)$.
\end{lem}
\proof 
We shall restrict the proof to the case of rightmost characteristics, the leftmost part being analogous. 

To begin the proof let us consider a point $t_0>0$ such that $\dot{x_\alpha}(t_0)=u(x_\alpha(t_0),t_0)$ and that the characteristic $x_\alpha$ branches at this point. The set of values of characteristics emanating from the branching point $(x_\alpha(t_0),t_0)$
\[
\{(z,t):z=x_\alpha(t), x_\alpha(0)=x_0, t_0\leq t\leq T\}
\] 
is compact and connected by Proposition \ref{kneser}. For any fixed $t_0\leq t\leq T$ one can take $y(t)$ as $\max$ of the elements of this set. By Proposition \ref{prop2.0.1}, $y(t)$ is also Lipschitz continuous. To prove the claim it is enough to show that $y(t)$ satisfies \eqref{2.0.2} at the points of differentiability. Indeed, by \cite[Lemma 3.1]{daf_jhde} we see that then $y(t)$ is $C^1$ regular.

Suppose, on the contrary, that $\dot{y}(t)\neq u(y(t))$ for some $t\in (t_0, T)$, which is a point of differentiability
of $y$. Without loss of generality, we may assume that 
\[
\dot{y}(t)\leq u(y(t))-\delta
\]  
for some $\delta>0$. Next, we notice that by the continuity of $u$, we can choose $\varepsilon>0$ in such a way that for every $(x,t)\in [y(t)-\varepsilon,y(t)+\varepsilon]\times [t-\varepsilon,t+\varepsilon]$ we have
\begin{equation}\label{2.0.5}
u(x,t)>\dot{y}(t)+\frac{\delta}{2},
\end{equation}
as well as $t+\varepsilon<T$, moreover for every $s\in [t-\varepsilon, t+\varepsilon]$
\begin{equation}\label{2.0.6}
|y(s)-y(t)-\dot{y}(t)(s-t)|<\frac{\delta(s-t)}{10}.
\end{equation}
By Proposition \ref{kneser} we can choose $\alpha_0$ such that $y(t)=x_{\alpha_0}(t)$. Then for $\tau$ small enough
\[
x_{\alpha_0}(t+\tau)=x_{\alpha_0}(t)+\int_t^{t+\tau}u(x_\alpha(s),s)ds\stackrel{\eqref{2.0.5}}{\geq} y(t)+\tau \left(\dot{y}(t)+\frac{\delta}{2}\right)\stackrel{\eqref{2.0.6}}{>}y(t+\tau),
\] 
contradiction.

\qed

Basing on the above lemma we define leftmost and rightmost characteristics. Finally, we can proceed with the proof of Lemma \ref{lemat2}.

\begin{figure}[htbp]
\begin{center}
\includegraphics[width=9cm]{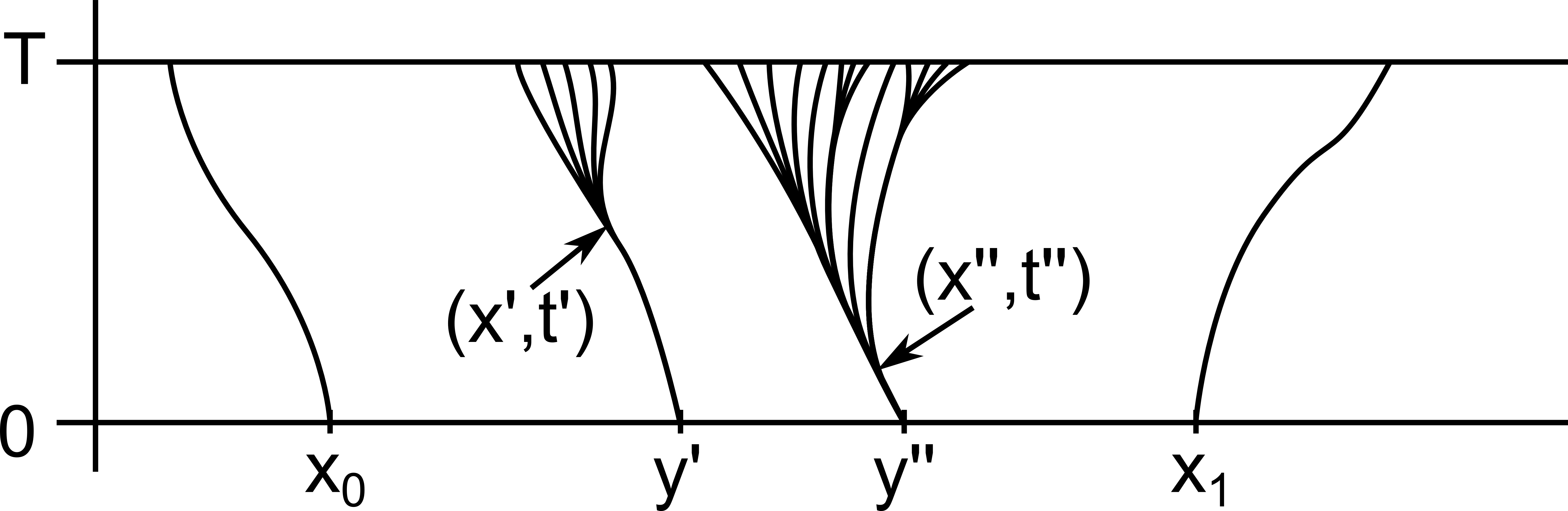}
\caption{Schematic presentation of branching sets related to $(x',t')$ and $(x'',t'')$ which are unions of graphs, until time $T$, of all characteristics emanating from $(x',t')$ and $(x'',t'')$, respectively.  Since $y', y''  \in I_T$, where $y',y''$ satisfy $x_{y'}(t')=x'$ and $x_{y''}(t'')=x''$,  these branching sets are disjoint.}
\label{Fig_branching}
\end{center}
\end{figure}

\vspace{0.3cm}
\textbf{Proof of Lemma \ref{lemat2}.}
Suppose that a characteristic $x(t)$ starting from $I_T$ branches for the first time at the point $(x(t_0),t_0)$. By Proposition \ref{kneser} and Lemma \ref{lemat3} the rightmost and leftmost characteristics emanating from the point $(x(t_0),t_0)$
together with the line $t=T$ bound a set of positive measure. We name such a set a branching set related to $(x(t_0),t_0)$. Consider now the interval $[x_0, x_1] \subset \mathbb{R}$ and all the characteristics emanating from $[x_0, x_1]\cap I_T$. First, we notice that by Corollary \ref{cor2.1} branching sets related to different points $(x',t'), (x'',t'')$ are disjoint, see Fig. \ref{Fig_branching}. Next, we claim that the set of points of first branching times $(x',t')$ is countable. Indeed, otherwise the measure of the branching sets related to all the branching points $(x',t')$ would be infinite, but this set is a subset of the set bounded by the interval $[x_0, x_1]$ from the bottom, by the curves $x_l[x_0]$ and $x_r[x_1]$-respectively the leftmost characteristic emanating from $x_0$ and the rightmost emanating from $x_1$, and the interval $[x_l[x_0](T),x_r[x_1](T)]$ from the top. The latter set is however of finite measure.   

\qed

In view of Corollary \ref{cor2.2} we define the following sets.
\begin{defi}\label{defi2.0.1}
By $I_t^{unique}$ we understand the full-measure subset of $I_t$ consisting of all points $x$ such that every characteristic of \eqref{hs} associated to the weak solution $u$ starting at $x$ stays unique up to time $t,\; t\in(0,\infty]$. 
\end{defi}

\mysection{Time-monotonicity of positive part of the energy and consequences}\label{section3}

We first define the positive part of the energy.
\begin{defi}\label{3.01}
Let $w_+(y,t):=\max \{w(y,t),0\}$ be the positive part of $w$. We define the positive part of the energy as 
\begin{equation}\label{3.0.1}
E^+_{[\eta,\xi]}(t):=\int_\eta^\xi w_+^2(y,t)dy.
\end{equation}
\end{defi}
Next, we recall that at $I_\infty$ the function $w_0(y)\geq 0$.
As is seen from the following proposition, the positive part of the energy defined above is nondecreasing.
\begin{prop}\label{prop3.0.1}
Let $u$ be a weak solution to \eqref{hs}, $w=u_x$ a.e. Take $a,b\in I_\infty^{unique}$ such that $a<b$. Then for every $t\in[0,\infty]$ we have
\begin{equation}\label{3.0.2}
E^+_{[x_a(t),x_b(t)]}(t)\geq E^+_{[a,b]}(0),
\end{equation}
where $x_a(t), x_b(t)$ are the unique characteristics emanating from $a$ and $b$, respectively.
\end{prop}
\proof
For every $\zeta\in I$ let $x_r[\zeta](t)$ denote the rightmost characteristic emanating from $\zeta$. Note that
for every $t > 0$ the mapping $M:\zeta\rightarrow x_r[\zeta](t)$ is monotone increasing. 
This mapping may be constant on some intervals (if characteristics from some interval meet before time $t>0$) or have jumps (case of branching before time $t>0$).

Take $\zeta\in I_\infty$. Next, take $\eta$ from a neighborhood of $\zeta$. By \eqref{omega} we have
\[
\omega^{x_\eta, x_\zeta}(t)\geq \frac{2\omega^{\eta, \zeta}(0)}{2+t\omega^{\eta, \zeta}(0)},
\]  
where we use the notation from the proof of Lemma \ref{lemat1}. Since $\zeta$ was chosen from $I_\infty$, if $\eta$ is close enough, then $\omega^{\eta, \zeta}(0)>0$, which in turn gives $\omega^{x_\eta, x_\zeta}(t)>0$. Consequently,
\begin{equation}\label{3.0.31}
I_\infty(t)\subset \{y: w(y,t)\;\mbox{exists and is non-negative}\;\},
\end{equation}
with 
\[
A(t):=\{x_r[a](t), a\in A\}\]
for $x_r[a](t)$ the rightmost characteristic starting from $a$, and
\begin{equation}\label{3.0.3}
w_+(y,t)\geq \liminf_{\eta\rightarrow\zeta}\omega^{x_\eta, x_\zeta}(t)\geq\liminf_{\eta\rightarrow\zeta}\frac{2\omega^{\eta, \zeta}(0)}{2+t\omega^{\eta, \zeta}(0)}=\frac{2w(\zeta,0)}{2+tw(\zeta,0)}.
\end{equation} 
On the other hand, in view of \eqref{3.0.31}, for $a,b\in I_\infty^{unique}$
\begin{equation}\label{3.0.33}
\int_{[x_a(t),x_b(t)]}w_+(y,t)^2dy\geq \int_{\left([a,b]\cap I_\infty\right)(t)}w_+(y,t)^2dy,
\end{equation}
where $x_a(t)$ is the characteristic starting from $a$.

Now observe that for $\zeta\in [a,b]\cap I_\infty$ equality $M(\zeta) = M(\eta)$ implies $\zeta= \eta$. Hence on the set
\begin{equation}\label{3.0.32}
M\left([a,b]\cap I_\infty\right)=\left([a,b]\cap I_\infty\right)(t)
\end{equation}
we can define a unique inverse mapping $M^{-1}$. This mapping can be prolonged to a right-continuous
generalized inverse of $M$ on $[a,b]$, which we call $W$. The definition of $W$, which can be taken for instance from \cite{embrechts}, reads
\[
W(y):=\inf\{x\in \R: M(x)\geq y\}, \;y\in [M(a), M(b)]. 
\]
Next, we are in a position to use the classical change of variable formula for the Lebesgue-Stieltjes integral, see
\cite[(1)]{falk_teschl}. We have
\begin{equation}\label{c_of_v} 
\int_a^b f(x)dM(x)=\int_{M(a)}^{M(b)}f(W(y))dy
\end{equation}
for any bounded Borel function $f:[a,b]\rightarrow \R$, any nondecreasing function $M:[a,b]\rightarrow \R$ and
its generalized inverse $W$. Choosing $f$ in \eqref{c_of_v} as 
\[
f(x):=\textbf{1}_{I_\infty}(x)\left(\frac{2w(x,0)}{2+tw(x,0)}\right)^2,
\]
which is nonnegative and bounded for fixed $t$, we arrive at
\begin{eqnarray}\label{3.0.4}
\int_{\left([a,b]\cap I_\infty\right)(t)}w_+(y,t)^2dy&\stackrel{\eqref{3.0.3}}{\geq}& \int_{[M(a),M(b)]}\textbf{1}_{I_\infty}(W(y))\left(\frac{2w(W(y),0)}{2+tw(W(y),0)}\right)^2dy\\ \nn
&=&\int_{[a,b]}f(x)dM(x)\geq \int_{[a,b]}f(x)M'(x)dx.
\end{eqnarray}
On the right-hand side of the above inequality we omitted integration over the singular part of the measure $dM$ due to positivity of $f$ and $M'(x)$ is computed a.e. Now,
observe that for every $\zeta\in I_\infty^{unique}$ we can estimate $M'(\zeta)$.
\begin{eqnarray}\label{M'}\nn
M'(\zeta)&=&\liminf_{\eta\rightarrow \zeta}\frac{x_\eta(t)-x_\zeta(t)}{\eta-\zeta}\stackrel{\eqref{wazne}}{=}\liminf_{\eta\rightarrow \zeta}e^{\int_0^t\omega^{x_\eta, x_\zeta}(s)ds}\stackrel{\eqref{omega}}{\geq} \liminf_{\eta\rightarrow \zeta}e^{\int_0^t\frac{2\omega^{\eta, \zeta}(0)}{2+s\omega^{\eta, \zeta}(0)}ds}
\\ 
&=&\liminf_{\eta\rightarrow \zeta}\frac{1}{4}\left(2+t\omega^{\eta, \zeta}(0)\right)^2 = \frac{1}{4}\left(2+tw(\zeta,0)\right)^2.
\end{eqnarray} 
Applying the above estimate in \eqref{3.0.4}, and using \eqref{3.0.33}
we arrive at
\[
\int_{[x_a(t),x_b(t)]}w_+(y,t)^2dy\geq \int_{[a,b]\cap I_\infty}\left(\frac{2w(x,0)}{2+tw(x,0)}\right)^2\frac{1}{4}\left(2+tw(x,0)\right)^2 dx=\int_{[a,b]}(w_0)_+^2(x)dx.
\]

\qed

We notice that choosing $a, b\notin I_\infty^{unique}$ we arrive at a similar conclusion as in the above proposition. Indeed, we have the following claims with less restrictive assumptions. 
\begin{cor}\label{cor3.1.1}
One can relax the assumptions of Proposition \ref{prop3.0.1}, assuming only that $a, b\in \R$. Then
\begin{equation}\label{3.1.0}
E^+_{[x_r[a](t),x_l[b](t)]}(t)\geq E^+_{[a,b]}(0),
\end{equation}
where $x_r, x_l$ stands for rightmost and leftmost characteristics. 
\end{cor}
\proof
Indeed, assume first $a\in I_\infty^{unique}, b\in\R$. Then there exists an increasing sequence $\left(b_{n}\right)_ {n\in \N}$ belonging to $I_\infty^{unique}$ such that $b_n\rightarrow \bar{b}$, where $\bar{b}:=\sup_{x<b, x\in I_\infty^{unique}}$. If $\bar{b}<b$ then we notice that
\begin{equation}\label{3.1.01}
w_0 \leq 0\;\mbox{a.e. on}\; [\bar{b}, b].
\end{equation}
Indeed, otherwise we would have a set $B\subset (\bar{b}, b]$ of positive measure such that $w_0(x)>0$ for all $x\in B$. Hence, $B\subset I_\infty$ and so it must have a nonempty intersection with $I_\infty^{unique}$, so there exists $x_0\in B\cap I_\infty^{unique}$ such that
$w_0(x_0)>0$, which contradicts the definition of $\bar{b}$.
We obtain
\[
\int_{[a,b]}(w_0)_+(x)^2dx\stackrel{\eqref{3.1.01}}{=}\lim_{n\rightarrow \infty}\int_{[a,b_n]}(w_0)_+(x)^2dx\stackrel{Prop.\ref{prop3.0.1}}{\leq}\lim_{n\rightarrow \infty}\int_{[x_a(t),x_{b_n}(t)]}w(y,t)_+^2dy
\]
\[
=\int_{[x_a(t),x_l[\bar{b}](t)]}w(y,t)_+^2dy\leq \int_{[x_a(t),x_l[b](t)]}w(y,t)_+^2dy.
\]
If both $a,b\in \R$, $a<b$ and $a,b\notin I_\infty^{unique}$, we find sequences $a_n, b_n\in I_\infty^{unique}$, respectively decreasing and increasing such that $a_n\rightarrow \bar{a}$, $b_n\rightarrow \bar{b}$, $\bar{a}:=\inf_{x>a, x\in I_\infty^{unique}}$. The same computation as above yields the claim. 

\qed

Moreover, we notice that repeating an adequate part of the proof of Proposition \ref{prop3.0.1} we arrive at the following fact.
\begin{cor}\label{cor3.1.2}
For $0<t<\infty$ and $\zeta\in I_t^{unique}$ it holds
\[
M'(\zeta)\geq\frac{1}{4}\left(2+tw(\zeta,0)\right)^2.
\]
\end{cor}
Finally, we are in a position to define $J\subset I$ of a full measure by
\[
J:=\{\zeta\in I:\;\mbox{ for a.e.}\; t\in [0,T_\zeta)\; (x_\zeta(t),t)\in \Gamma\},
\] 
where $x_\zeta(t)$ is a characteristic, $T_\zeta$ was defined in the Definition \ref{I_s},  and 
\[
\Gamma:=\{(x,t):w=\partial_xu(x,t) \;\mbox{exists}\;\}.
\]
By the definition, for $t\in (0, \infty]$ $J_t^{unique}:=I_t^{unique}\cap J$. 

The next lemma is crucial in our proof of Theorem \ref{Tw.1}. It allows us to control the difference quotients $\omega$ on a subset of $J_t^{unique}$ of full measure.  
\begin{lem}\label{lem3.1.1}
Let $u$ be any weak solution to \eqref{hs}. Fix $0<t<\infty$ and $0<\tau<t$. For almost every $\zeta\in J_t^{unique}$
there exist $M > 0$ and $\varepsilon>0$ such that for every $\eta \in (\zeta,\zeta+\eps]$  we have $\omega^{x_\zeta,x_l[\eta]}(s)\le M$ for every $0\leq s\leq \tau$. 
\end{lem}
\proof
We denote by $J_t^{bad}$ the following set.
\[
J_t^{bad}:=\{\zeta\in J_t^{unique}:\forall \varepsilon>0\; \forall M>0 \;\exists \eta:\eta-\zeta<\varepsilon \;\exists s\in [0, \tau]\;\mbox{such that}\; \omega^{x_\zeta, x_l[\eta]}(s)>M \}.
\]
In order to prove Lemma \ref{lem3.1.1} it is enough to show that the measure of $J_t^{bad}$ is zero. To this end
fix $M>0$ and for $\zeta\in J_t^{bad}$ denote
\[
\Pi_\zeta^{M,\delta}:=[\zeta, \eta],
\]
where $\eta$ is a point that satisfies
\begin{equation}\label{3.0.45}
\eta\in (\zeta, \zeta+\delta)\;\mbox{and there exists}\;0\leq s\leq \tau\;\mbox{such that}\;\omega^{x_\zeta, x_l[\eta]}(s)>M.
\end{equation}
We observe that for fixed $M>0$ 
\[
{\cal E}^M:=\{\Pi_\zeta^{M,\delta}, \delta>0, \zeta\in J_t^{bad}\}
\]
is a covering of $J_t^{bad}$. Moreover, any point in $J_t^{bad}$ is contained in an element of ${\cal E}^M$ of arbitrarily small length. Indeed, by the definition of $J_t^{bad}$ one sees that given $\zeta\in J_t^{bad}$ for any small $\delta>0$ there exists $\eta \in(\zeta, \zeta+\delta)$ for which $\omega^{x_\zeta, x_l[\eta]}(s)>M$. So, ${\cal E}^M$ is a Vitali covering of $J_t^{bad}$. By the Vitali theorem, we obtain at most countable subfamily ${\cal F}^M\subset {\cal E}^M$ of pairwise disjoint closed intervals such that 
\[
J_t^{bad}\subset \bigcup {\cal F}^M
\]
holds up to a set of measure $0$. Denote 
\[
{\cal F}^M:=\{[\zeta_i, \eta_i]\}, i\in \N, \zeta_i\in J_t^{bad}, \eta_i \;\mbox{satisfy}\;\eqref{3.0.45}.
\]
Then, for any $i\in \N$ there exists $s_i\in [0, \tau]$ such that $\omega^{x_{\zeta_i}, x_l[\eta_i]}(s)>M$. This leads us to 
\begin{equation}\label{3.0.451}
\int_{x_{\zeta_i}(s_i)}^{x_l[\eta_i](s_i)} w_+^2(y,s_i)dy\geq (x_l[\eta_i](s_i)-x_{\zeta_i}(s_i))M^2.
\end{equation}
Indeed, \eqref{3.0.451} holds since by the Schwarz inequality and in view of the obvious inequality $w_+\geq w$ 
\begin{eqnarray*}
\int_{x_{\zeta_i}(s_i)}^{x_l[\eta_i](s_i)} w_+^2(y,s_i)dy&\geq&  \frac{1}{x_l[\eta_i](s_i)-x_{\zeta_i}(s_i)}\left(\int_{x_{\zeta_i}(s_i)}^{x_l[\eta_i](s_i)} w_+(y,s_i)dy\right)^2 \\
&\geq &\frac{\left(u(x_l[\eta_i](s_i))-u(x_{\zeta_i}(s_i))\right)^2}{x_l[\eta_i](s_i)-x_{\zeta_i}(s_i)}>(x_l[\eta_i](s_i)-x_{\zeta_i}(s_i))M^2.
\end{eqnarray*}
In view of Proposition \ref{prop3.0.1} and Corollary \ref{cor3.1.1}, \eqref{3.0.451} yields 
\[
\int_{x_{\zeta_i}(\tau)}^{x_l[\eta_i](\tau)} w_+^2(y,\tau)dy\geq \int_{x_{\zeta_i}(s_i)}^{x_l[\eta_i](s_i)} w_+^2(y,s_i)dy \geq (x_l[\eta_i](s_i)-x_{\zeta_i}(s_i))M^2.
\]
Summing over $i\in \N$ we obtain
\begin{equation}\label{3.0.47}
\int_{\R}w^2(y,\tau)dy\geq \int_{\R}w_+^2(y,\tau)\geq \sum_{i=1}^\infty  \int_{x_{\zeta_i}(\tau)}^{x_l[\eta_i](\tau)} w_+^2(y,\tau)dy   \geq M^2\sum_{i=1}^\infty (x_l[\eta_i](s_i)-x_{\zeta_i}(s_i)).
\end{equation}
We observe the following estimate
\begin{equation}\label{3.0.475}
(x_l[\eta_i](s_i)-x_{\zeta_i}(s_i))\geq \frac{(\eta_i-\zeta_i)\left(1-\frac{\tau}{t}\right)^2}{2}\;.
\end{equation}
Indeed, for $h(s)=x_l[\eta_i](s)-x_{\zeta_i}(t)$, as long as $h>0$, \eqref{wazne} is satisfied. Hence
\[
h(t)\geq (\eta_i-\zeta_i)e^{\int_0^t\omega^{x_{\zeta_i}, x_l[\eta_i]}(s)ds}.
\]
On the other hand, one can estimate 
\[
\omega^{x_{\zeta_i}, x_l[\eta_i]}(s)\geq \frac{2\omega^{\zeta_i, \eta_i}(0)}{2+s\omega^{\zeta_i,\eta_i}(0)}
\]
the same way as in \eqref{omega}. Consequently,
\begin{equation}\label{raz_dwa}
h(s_i)\geq \frac{\eta_i-\zeta_i}{4}\left(2+s_i\omega^{\zeta_i,\eta_i}(0)\right)^2\geq \frac{\eta_i-\zeta_i}{2}\left(1-\frac{\tau}{t}\right)^2,
\end{equation}
where in the last inequality we made use of the inequalities $s_i<\tau$ and $\omega^{\zeta_i, \eta_i}(0)>-\frac{2}{t}$, the latter holds on $J_t$.
Plugging \eqref{3.0.475} in \eqref{3.0.47} we arrive at
\[
\int_{\R}w^2(y,\tau)dy\geq M^2\sum_{i=1}^\infty \frac{(\eta_i-\zeta_i)\left(1-\frac{\tau}{t}\right)^2}{2}\geq \frac{M^2}{2}\left(1-\frac{\tau}{t}\right)^2|J_t^{bad}|.
\]
But the last inequality means
\[
|J_t^{bad}|\leq \frac{2\int_{\R}w^2(y,\tau)dy}{M^2\left(1-\frac{\tau}{t}\right)^2}\;,
\]
so enlarging $M$, we see that $|J_t^{bad}|=0$.

\qed

\mysection{Maximal dissipation selects the unique solution}\label{section4}

The present section consists of two subsections. In the first one we prove an averaging lemma, which will be used in
the second one in order to prove Proposition \ref{prop2.3}.

\subsection{Averaging lemma}

We prove the following proposition.
\begin{prop}\label{prop4.1.1}
For any $g\in L^1(\R)$ the following formula holds
\begin{equation}\label{4.1.1}
\lim_{\varepsilon\rightarrow 0} \int_{\R}\frac{1}{\varepsilon}\int_x^{x+\varepsilon}|g(y)-g(x)|dydx=0.
\end{equation}
\end{prop}
\proof
By the Fubini theorem we obtain
\begin{eqnarray*}
\int_{\R}\frac{1}{\varepsilon}\int_x^{x+\varepsilon}|g(y)-g(x)|dydx&=&\int_{\R}\frac{1}{\varepsilon}\int_0^\varepsilon|g(x+y)-g(x)|dydx\\ \nn
&=&\int_0^\varepsilon\frac{1}{\varepsilon}\left(\int_{\R}|g(x+y)-g(x)|dx\right)dy\\ \nn
&\leq& \sup_{y\in [0,\varepsilon] }\left(\int_{\R}|g(x+y)-g(x)|dx\right).\nn
\end{eqnarray*}
By the continuity of the translation in $L^1$ we infer 
\[
\lim_{\varepsilon\rightarrow 0}\sup_{y\in [0,\varepsilon] }\left(\int_{\R}|g(x+y)-g(x)|dx\right)=0
\]
which yields the claim.

\qed

\subsection{Averaging over characteristics}


Following Lemma \ref{lem3.1.1} and definition of $I_t$, we can represent, up to a set of measure zero, $J_T^{unique}$ as a countable union of sets with $\omega$ bounded on $[0,\tau]$ (this property will be crucial in our proof and is the fundament of the decomposition of $J_T^{unique}$ which we introduce) and with $\omega(0)$ close to $w(0)$. More precisely, there exists a set $Z$ of measure $0$ such that
\[ 
J_T^{unique} = \left(\bigcup_{N=1}^{\infty} J_{T}^{unique,N}\right) \cup Z, 
\]
where 
\[ J_T^{unique, N} := \{\zeta \in J_T^{unique}: \omega^{x_\zeta, x_l[\zeta+\eps]}(s) \le N \mbox{ and } \omega^{\zeta, \zeta+\eps}(0) \ge -2 \slash T   \mbox{ for } \eps\le \frac {1}{N} \mbox{ and } 0\le s \le \tau\}\],  
and we used Lemma \ref{lem3.1.1} as well as the fact that $w_0(\zeta) > -2 \slash T$.

\begin{prop}\label{Cor_constants}
Let $0 \le t \le \tau <T$. If $\zeta \in J_T^{unique,N}$ then for $M(\zeta):=x_\zeta(t)$ and $\eps\le 1\slash N$ we have:
\begin{enumerate}[i)]
\item $M'(\zeta) \ge c(\tau)$, 
\item $M(\zeta+\eps) - M(\zeta) \le \eps C_N(\tau)$, 
\item $M(\zeta+\eps) - M(\zeta) \ge \eps {c}(\tau)$,
\end{enumerate}
where $c(\tau) := \frac 1 2 (1  -\tau\slash T)^2 $ and $C_N(\tau):=e^{N\tau}.$
\end{prop}
\proof
The proof of parts (i) and (iii) consists of repeating the argument in \eqref{M'} in the context of the present proposition. Since $\tau<T$ we arrive at the desired claim. As a consequence of \eqref{wazne} we obtain (ii). 

\qed

Below we formulate and prove a result which is a slight extension of the Riesz lemma on choosing an a.e. convergent subsequence from a sequence convergent in $L^p$.
\begin{prop}\label{riesz}
Let $(X,\mu)$ be a measure space and let $D_1 \subset D_2 \subset \dots$ be an increasing family of subsets of $X$. Let $D\subset X$ satisfy $\mu(D \backslash  \bigcup_{n=1}^{\infty} D_n) = 0$ and $\mu(\bigcup_{n=1}^{\infty} D_n \backslash D) = 0$. Consider a family of functions $d_\eps(\zeta):X\rightarrow \R$, $\eps\in(0, \eps_0)$, such that $d_\eps(\zeta)\stackrel{\eps\rightarrow 0}{\longrightarrow} 0$ in $L^1(D_n)$ for $n=1,2,\dots$. Then there exists a subsequence $\eps_k$ tending to zero such that
\[
d_{\eps_k}(\zeta)\rightarrow 0 \;\;a.e. \;\mbox{in}\; D. 
\]  
\end{prop}
\proof
We use the diagonal argument. Namely, convergence in $L^1$ implies convergence almost everywhere on a subsequence. Hence, there exists a convergent to $0$ and decreasing sequence 
$\eps^1_1, \eps^1_2, \dots$  such that $d_{\eps^1_k} \stackrel{k \rightarrow 0}{\longrightarrow} 0$ a.e. on $D_1$. Define inductively sequence
$\eps^n_1,\eps^n_2, \dots$
as a subsequence of $\eps^{n-1}_2,\eps^{n-1}_3, \dots $ satisfying $d_{\eps^n_k} \stackrel{k \rightarrow 0}{\longrightarrow} 0$ a.e. on $D_n$. 
Finally, take $\eps_k := \eps^k_1$ for $k = 1,2,\dots$. Since $\{\eps_k\}_{k=n}^{\infty}$ is a subsequence of $\{\eps^n_j\}$ for every $n$, we obtain
\[
d_{\eps_k}(\zeta)\rightarrow 0 \;\;a.e. \;\mbox{in}\; D_n
\]  
for every $n$. Since almost every $\zeta \in D$ belongs in fact to some $D_n$, we conclude.
\qed

Let us now state and prove a crucial lemma on averaging the energy over characteristics. 
\begin{lem}
\label{Lem_42}
Let $u$ be a weak solution to \eqref{hs}, $w=u_x$ a.e. Let $[0,\tau] \subset [0,T)$. Then for almost every $\zeta \in J_T^{unique}$ and every $0\le \sigma \le \tau$ we have
\begin{equation}
\label{Eq_conveps}
\limsup_{k \to \infty}  \left| \int_{\sigma}^{\tau} \left(  \frac {1}{x_l[\zeta+\eps_k](t) - x_\zeta(t)} \int_{x_\zeta(t)}^{x_l[\zeta+\eps_k](t)} w^2(y,t)dy - w^2(x_\zeta(t),t) \right)    dt\right| = 0,
\end{equation}
where $\eps_k$ is some sequence convergent monotonically to $0$ and $x_l[\zeta+\eps_k](t)$ is the leftmost characteristic emanating from $\zeta + \eps_k$, associated to $u$.
\end{lem}
\proof
It is enough to show that for every $N  \in \mathbb{N}$ we have
\begin{equation}
\label{Eq_toshow}
\lim_{\eps \to 0^+}  \int_{J_T^{unique,N}}  \int_{\sigma}^{\tau} \left( \frac {1}{x_l[\zeta+\eps](t) - x_\zeta(t)} \int_0^{x_l[\zeta+\eps](t)-x_\zeta(t)}|w^2(x_\zeta(t) + y,t) - w^2(x_\zeta(t),t)| dy \right) dt d\zeta = 0.
\end{equation}

Indeed, one applies Proposition \ref{riesz} with $X=\R$, $\mu$ being a Lebesgue measure, $D=J_T^{unique}, D_N=J_T^{unique,N}$ and 
\[
d_\eps(w,\zeta):=\int_{\sigma}^{\tau} \left( \frac {1}{x_l[\zeta+\eps](t) - x_\zeta(t)} \int_0^{x_l[\zeta+\eps](t)-x_\zeta(t)}|w^2(x_\zeta(t) + y,t) - w^2(x_\zeta(t),t)| dy \right) dt.
\]
It remains to show \eqref{Eq_toshow}. To this end, first observe that for $\eps \le 1 \slash N$
\begin{eqnarray*}
&&\int_{J_T^{unique, N}} \left( \frac {1}{x_l[\zeta+\eps](t)-x_\zeta(t)} \int_0^{x_l[\zeta+\eps](t)-x_{\zeta}(t)}|w^2(x_\zeta(t) + y,t) - w^2(x_\zeta(t),t)| dy \right) d\zeta \\
&=& \int_{J_T^{unique, N}}  \left( \frac {1}{M(\zeta+\eps)-M(\zeta)} \int_0^{M(\zeta+\eps)-M(\zeta)}|w^2(M(\zeta) + y,t) - w^2(M(\zeta),t)| dy \right)d\zeta \\ &\stackrel{\eqref{Cor_constants}\; ii), iii)}{\le}& 
\int_{J_T^{unique, N}}  \left( \frac {1}{c(\tau)\eps} \int_0^{C_N(\tau)\eps}|w^2(M(\zeta) + y,t) - w^2(M(\zeta),t)| dy \right)d\zeta \\&=& 
\int_{J_T^{unique, N}} g^{\eps}(M(\zeta))d\zeta,
\end{eqnarray*}
where 
\[
g^{\eps}(z):= \textbf{1}_{M(J_T^{unique,N})}(z)\frac {1}{c(\tau)\eps} \int_0^{C_N(\tau)\eps} |w^2(z+y,t) - w^2(z,t)|dy. 
\]
Now, observe that for fixed $\eps \le 1 \slash N$ function $g^\eps$ is nonnegative, bounded and Borel measurable.
Using \cite[(6)]{falk_teschl}, we obtain
\[
\int_{M\left(J_T^{unique,N}\right)} g^\eps(z)dz = \int_{J_T^{unique,N}} g^\eps(M(\zeta))dM(\zeta) 
\]
Next note that neglecting the singular part of $dM$, using nonnegativity of $g^\eps$ we continue
\[  
\geq \int_{J_T^{unique,N}} g^\eps(M(\zeta))M'(\zeta)d\zeta \stackrel{\eqref{Cor_constants}\; i)}{\ge} c(\tau) \int_{J_T^{unique,N}} g^\eps(M(\zeta))d\zeta.
\]
Consequently, 
\begin{eqnarray*}
S^{\eps}(t) &:=& \int_{J_T^{unique,N}} g^\eps(M(\zeta))d\zeta\\ &\le&
\frac {1}{c(\tau)} \int_{M\left(J_T^{unique,N}\right)} g^\eps (z)dz\\ &=&
\frac {1}{c(\tau)}\int_{M\left(J_T^{unique,N}\right)}  \frac {1}{c(\tau)\eps} \int_0^{C_N(\tau)\eps} |w^2(z+y,t) - w^2(z,t)|dydz\\ &=&
\frac {C_N(\tau)}{c^2(\tau)}\int_{M\left(J_T^{unique,N}\right)}  \frac {1}{C_N(\tau)\eps} \int_0^{C_N(\tau)\eps} |w^2(z+y,t) - w^2(z,t)|dydz\\ &\le& 
\frac {C_N(\tau)}{c^2(\tau)}\int_{\mathbb{R}}  \frac {1}{C_N(\tau)\eps} \int_0^{C_N(\tau)\eps} |w^2(z+y,t) - w^2(z,t)|dydz,
\end{eqnarray*}
where constants $c(\tau), C_N(\tau)$ are defined in Proposition \ref{Cor_constants}. Using Proposition \ref{prop4.1.1} we see that 
\[
S^\eps(t) \to 0 
\] 
as $\eps \to 0$ for almost every $t\in [0,\tau]$. Note also that, by the Fubini theorem, 
\[ S^\eps(t) \le 2 \frac {C_N(\tau)}{c^2(\tau)} \sup_{t \in [0,T)} E(t). \] Using the Fubini theorem once again as well as the Lebesgue dominated convergence theorem, we obtain

\begin{eqnarray*}
\lim_{\eps \to 0^+}  \int_{J_T^{unique,N}}  \int_{\sigma}^{\tau} \left( \frac {1}{x_l[\zeta+\eps](t) - x_\zeta(t)} \int_0^{x_l[\zeta+\eps](t)-x_\zeta(t)}|w^2(x_\zeta(t) + y,t) - w^2(x_\zeta(t),t)| dy \right) dt d\zeta \\
=\lim_{\eps \to 0^+} \int_\sigma^\tau S^\eps(t)dt = 0.
\end{eqnarray*}

This proves \eqref{Eq_toshow}.

\qed

Thus, we have completed all the preparatory steps and now we proceed with the proof of the main result. 

\vspace{0.3cm}
\textbf{Proof of Theorem \ref{Tw.1}}

As it was explained in Section \ref{section1.5} in order to prove Theorem \ref{Tw.1}, it is enough to prove Proposition \ref{prop2.3}. It implies Propositions \ref{prop2.1} and \ref{prop2.2}, which in turn yield the main theorem. 

Let the sequence $\eps_k$ be obtained as in Lemma \ref{Lem_42}. 
Observe that $J_T^{unique}= \bigcup_{K=1}^{\infty} J_{T + 1 \slash K}^{unique}$ (up to a set of measure $0$).  Now fix $K \in \mathbb{N}$. For $\zeta \in J_{T+1\slash K}^{unique}$ and $t\leq T$, by \eqref{omegaintegral}, we have
\[ \frac {d}{dt}{\omega}^{x_\zeta, x_l[\zeta+\eps_k]} (t) = -\omega^{x_\zeta, x_l[\zeta+\eps_k]}(t)^2 + \frac {1}{2 (x_l[\zeta+\eps_k](t) - x_\zeta(t))}  \int_{x_\zeta(t)}^{x_l[\zeta+\eps_k](t)} w^2(y,t)dy. \]


Hence,  for $0 \le \sigma \le t \le \tau < T + 1\slash K$ we obtain
\begin{equation}
\label{eq_limitomega}
{\omega}^{x_\zeta, x_l[\zeta+\eps_k]} (\tau) - {\omega}^{x_\zeta, x_l[\zeta+\eps_k]} (\sigma)  = - \int_\sigma^\tau \omega^{x_\zeta, x_l[\zeta+\eps_k]}(t)^2 dt + \frac {1}{2}  \int_\sigma^\tau \frac{1}{x_l[\zeta+\eps_k](t) - x_\zeta(t)}  \int_{x_\zeta(t)}^{x_l[\zeta+\eps_k](t)} w^2(y,t)dy dt.
\end{equation}

Passing to the limit $k\to\infty$ in the left-hand side of \eqref{eq_limitomega} we obtain 
\begin{equation}
\label{Eq_LHS}
\lim_{k \to \infty} \left({\omega}^{x_\zeta, x_l[\zeta+\eps_k]} (\tau) - {\omega}^{x_\zeta, x_l[\zeta+\eps_k]} (\sigma)\right)  =  w(x_\zeta(\tau),\tau) - w(x_\zeta(\sigma),\sigma) 
\end{equation}
for almost every $\sigma,\tau$ (see the definition of $\Gamma$ above Lemma \ref{lem3.1.1}). To pass to the limit in the right-hand side, we first observe that for almost every $\zeta \in J_{T+ 1 \slash K}^{unique}$ there exists $N_\zeta$ such that $\zeta \in J_{T+ 1 \slash K}^{unique,N_\zeta}$. Hence, for $k$ large enough 
\[
\omega^{x_\zeta, x_l[\zeta+\eps_k]}(s) \le N_\zeta
\] 
for $s \in [0,\tau]$. On the other hand, by \eqref{omega}
\[
\omega^{x_\zeta, x_l[\zeta+\eps_k]}(s) \ge \frac {-2\slash (T + 1 \slash K)}{1-(\tau\slash (T+ 1 \slash K))} 
\]
for $s\in [0,\tau]$ and  $k$ large enough.
Hence, $ \omega^{x_\zeta, x_l[\zeta+\eps_k]}(t)^2$ is bounded on $[0,\tau]$ and using the Lebesgue dominated convergence theorem we obtain 
\begin{equation}
\lim_{k\to \infty} \int_\sigma^\tau \omega^{x_\zeta, x_l[\zeta+\eps_k]}(t)^2 dt = \int_\sigma^\tau w^2(x_\zeta(t),t)dt.
\end{equation}
Finally, by Lemma \ref{Lem_42} 
\begin{equation}
\label{Eq_final}
\lim_{k\to \infty} \frac 1 2 \int_\sigma^\tau \frac{1}{x_l[\zeta+\eps_k](t) - x_\zeta(t)}  \int_{x_\zeta(t)}^{x_l[\zeta+\eps_k](t)} w^2(y,t)dy dt = \frac 1 2 \int_\sigma^\tau w^2(x_\zeta(t),t) dt
\end{equation}
Combining  \eqref{eq_limitomega}-\eqref{Eq_final} and summing over $K \in \mathbb{N}$, for almost every $\zeta \in J_T^{unique}$  we have  
\[ 
w(x_\zeta(\tau),\tau) - w(x_\zeta(\sigma),\sigma) =  -\frac{1}{2} \int_{\sigma}^{\tau} w^2(x_\zeta(t),t) dt 
\]
for almost every $0 \le \sigma \le \tau \le T$ (more precisely, for those $0 \le \sigma \le \tau \le T$ for which $w(x_{\zeta}(\sigma),\sigma)$ and $w(x_{\zeta}(\tau),\tau)$ exist).
Solving the above differential equation for a.e. $\zeta \in J_T^{unique}$, we obtain that
\begin{equation}\label{5.0.1}
w(x_\zeta(t),t) = \frac {2w_0(\zeta)}{2+tw_0(\zeta)}  
\end{equation}
for those $t \in [0,T]$ for which $w(x_{\zeta}(t),t)$ exists. In particular, this holds for $t=T$ and almost every $\zeta \in J_T^{unique}$.
Since $x_a(T)\leq x_r[a](T)$ and $x_b(T)\geq x_l[b](T)$ by the definition of rightmost and leftmost characteristics we arrive at
\begin{eqnarray*}
\int_{x_a(T),x_b(T)}w^2(y,T)dy &\ge& \int_{[x_r[a](T),x_l[b](T)]}w^2(y,T)dy \ge  \int_{M([a,b] \cap J_T^{unique})} w^2(y,T)dy\\ &\stackrel{\eqref{5.0.1}}{\ge}& \int_{M([a,b] \cap J_T^{unique})} \left( \frac {2 w_0(W(y))}{2+Tw_0(W(y))}\right)^2 dy \\ 
 &\stackrel{\eqref{c_of_v}}{=}& \int_{[a,b] \cap J_{T} ^{unique}}   \left( \frac {2 w_0(\zeta)}{2+Tw_0(\zeta)}\right)^2 dM(\zeta) \\
 &\stackrel{Cor.\ref{cor3.1.2}}{\ge}& \int_{[a,b] \cap J_{T} ^{unique}}   \left( \frac {2 w_0(\zeta)}{2+Tw_0(\zeta)}\right)^2  \times \frac {1}{4} [2+Tw_0(\zeta)]^2 d\zeta \\
  &=& \int_{[a,b] \cap J_T^{unique}} w_0(\zeta)^2 d\zeta.
\end{eqnarray*}
The claim of Proposition \ref{prop2.3} follows in view of the fact that $I_T=J_T^{unique}$ up to a set of measure zero. This, in turn, implies Theorem \ref{Tw.1}.

\qed

\noindent
{\bf Acknowledgement.} T.C. was partially supported by the National Centre of Science (NCN) under grant 2013/09/D/ST1/03687.

\end{document}